\documentclass[usebibtex]{gammas}

\newcount\Comments  
\newcommand{\kibitz}[2]{\ifnum\Comments=1\textcolor{#1}{#2}\fi}

\makeatletter
\renewcommand*\env@matrix[1][*\c@MaxMatrixCols c]{%
  \hskip -\arraycolsep
  \let\@ifnextchar\new@ifnextchar
  \array{#1}}
\makeatother
\DeclareMathAlphabet{\mathcal}{OMS}{cmsy}{m}{n}
\newcommand{\M}{\mathcal{M}}
\newcommand{\trace}{\text{tr}}
\newcommand{\N}{\mathbb{N}}
\newcommand{\Z}{\mathbb{Z}} 
\newcommand{\R}{\mathbb{R}} 

\newcommand{\Skew}{\operatorname{skew}} 
\newcommand{\dist}{\operatorname{dist}}

\newtheorem{proposition}{Proposition}
\newtheorem{remark}{Remark}
\newtheorem{corollary}{Corollary}

\gammtitle{On the Injectivity Radius of the Stiefel Manifold: Numerical investigations and an explicit construction of a cut point at short distance}

\gammauthora{Jakob Stoye\inst{a}\corauth}
\gammauthorb{Ralf Zimmermann\inst{b}\corauthb}

\gammauthoraorcid{0009-0003-9119-0013}
\gammauthorborcid{0000-0003-1692-3996}

\gammauthorhead{J.~Stoye, R.~Zimmermann}

\gammaddressa{Technische Universität Braunschweig, Braunschweig, Germany}	%
\gammaddressb{University of Southern Denmark, Department of Mathematics and Computer Science,
Odense, Denmark} 

\gammcorrespondence{\href{jakob.stoye@tu-braunschweig.de}{jakob.stoye@tu-braunschweig.de}}
\gammcorrespondenceb{\href{zimmermann@imadasdu.dk}{zimmermann@imadasdu.dk}, \href{https://portal.findresearcher.sdu.dk/en/persons/zimmermann}{https://portal.findresearcher.sdu.dk/en/persons/zimmermann}}

\usepackage{blindtext}

\newtheorem{theorem}{Theorem}
\newtheorem{lemma}[theorem]{Lemma}

\theoremstyle{definition}
\newtheorem{definition}[theorem]{Definition}

\begin{document}
\maketitle
%

\begin{gammabstract}
    Arguably, geodesics are the most important geometric objects on a differentiable manifold.
    They describe candidates for shortest paths and are guaranteed to be unique shortest paths when the starting velocity stays within the so-called injectivity radius of the manifold. 
    In this work, we investigate the injectivity radius of the Stiefel manifold under the canonical metric.
    The Stiefel manifold $St(n,p)$ is the set of rectangular matrices of dimension $n$-by-$p$ with orthonormal columns, sometimes also called the space of orthonormal $p$-frames in $\R^n$.
    Using a standard curvature argument, Rentmeesters \cite{Quentin} has shown that the injectivity radius of the Stiefel manifold is bounded by $\sqrt{\frac{4}{5}}\pi\approx 0.8944\pi$. It is an open question, whether this bound is sharp. With the definition of the injectivity radius via cut points of geodesics, we gain access to the information of the injectivity radius by investigating geodesics.
    More precisely, we consider the behavior of special variations of geodesics, called Jacobi fields.
    By doing so, we are able to present an explicit example of a cut point
    on a low-dimensional Stiefel manifold at a distance of ca. $0.9133\pi$. The precise value is given by the first positive root of $\left(\frac{t}{\sqrt{2}}\cos\left(\frac{t}{\sqrt{2}}\right) + \sin\left(\frac{t}{\sqrt{2}}\right)\right)$.
    In addition to the theoretical analysis, we investigate the question of the sharpness of the bound for the injectivity radius by means of numerical experiments. 
\end{gammabstract}

\begin{gammkeywords}
	Stiefel manifold, Injectivity radius, Riemannian Computing, Canonical metric, Cut points, Jacobi fields\\
\end{gammkeywords}

\section{Introduction}
The Riemannian manifold defined by the set of rectangular matrices with orthonormal columns is called the Stiefel manifold $St(n,p)$. Stiefel manifolds feature in a large variety of application problems, ranging from optimization
\cite{book_AMS, boumal2023, sato2021} over numerical methods for differential equations \cite{BennerGugercinWillcox2015, Celledoni_2020, HueperRollingStiefel2008, Zim21}
to applications in statistics and data science \cite{Turaga_2008, Chakraborty2018, FLETCHER2020}. 
On a manifold, the selected Riemannian metric determines how lengths and angles are measured, and thus how geodesics are defined. Geodesics give rise to a special set of local coordinate charts, the so-called Riemannian normal coordinates. These are the Riemannian exponential map and the Riemannian logarithm map.

In this work, we will consider the Stiefel manifold under the canonical metric. 
In this case, the geodesics are known in closed form \cite{EdelmanAriasSmith1998}.
For a starting point $U\in St(n,p)$ and a normalized starting velocity $\Delta\in T_U St(n,p)$ from the tangent space, the corresponding geodesic is given by the Stiefel exponential $\text{Exp}_U(t\Delta)$ at $U$. 
Geodesics are candidates for shortest paths and are unique shortest paths when the starting velocity stays within the so-called injectivity radius. The same condition ensures that the Stiefel exponential at $U$ and thus the Riemannian normal coordinates at that location are invertible. In this case, we are able to calculate the shortest path between two given points. Solving the geodesic endpoint problem is important, e.g., for interpolation tasks and for computing Riemannian centers of mass. For the injectivity radius on the Stiefel manifold a theoretical bound is given in \cite{Quentin}, 
\[
i(St(n,p)) \ge \sqrt{\frac{4}{5}}\pi.
\]
This bound stems from an extreme-case estimate of the sectional curvature on Stiefel that has been confirmed in \cite{ZimmermannStoye2024}.
One aim of this work is to investigate whether the bound on the injectivity radius is sharp, i.e., whether a geodesic can be found whose cut point is at $t=\sqrt{\frac{4}{5}}\pi$. For this purpose, we conduct numerical experiments with random geodesics of different lengths and investigate whether there are shorter geodesics to its start and end points. If there are shorter geodesics, the examined geodesic is no longer minimizing and therefore the cut point of the geodesic has already been reached. Hence, the injectivity radius must be smaller than the length of the geodesic whose cut point has already been reached. In the experiments, however, we are not able to reach the bound to the injectivity radius.\\
Furthermore, we construct an explicit example of a cut point on the Stiefel manifold $St(4,2)$. Here we consider a geodesic with velocity from a tangent plane section of maximal sectional curvature. For the construction, we calculate all linearly independent Jacobi fields along the geodesic under consideration and use this to compute the first conjugate point, which simultaneously describes the cut point of the geodesic. The distance of the cut point, i.e., the cut time, coincides with the upper bound on the injectivity radius observed in the numerical experiments. To the best of our knowledge, this is the first explicit presentation of a cut point on the Stiefel manifold in literature. 
This contributes to a better understanding of the geometry of the Stiefel manifold. The geometry of the manifold plays an important role in the efficient design of optimization algorithms and helps on the way to proving the injectivity radius.

\paragraph{Organization}
In \Cref{chap:background}, we recap general concepts on differentiable manifolds such as tangent spaces, Riemannian metrics, geodesics, Riemannian exponential and the injectivity radius as well as the geometry of quotient spaces. 
Then we recount these concepts for the concrete case of the Stiefel manifolds. 
In \Cref{sect:invoninj} we discuss the injectivity radius in detail. For this purpose, we introduce the concepts of curvature, Jacobi fields, conjugate points and cut points and review the bound for the injectivity radius of the Stiefel manifold from \cite{Quentin}. 
In \Cref{sec:numexpInj}, we conduct numerical experiments on whether the bound on the injectivity radius of the Stiefel manifold is sharp.
In \Cref{sect:explcutpoint}, we give an explicit example of a cut point on the Stiefel manifold $St(4,2)$. 
The same construction can be embedded in any Stiefel manifold of dimension $p\geq 2, n\geq p+2$. \Cref{sect:Summary} concludes the paper.

\paragraph{Notation}
For the reader's convenience, we list the main acronyms and variables.

\begin{small}
	\begin{tabular}{@{}ll@{}}
		Symbol &  meaning \\
		\midrule
		$I,I_n$   & identity matrix, provided with a dimensional index if required \\
		$\langle \cdot,\cdot\rangle,\,\|\cdot\|$ & canonical metric and the associated norm\\
		$O(n)$    & orthogonal group $O(n) = \{Q\in \R^{n\times n}\mid Q^TQ = I_n\}$\\
		$SO(n)$    & special orthogonal group $SO(n) = \{Q\in O(n)\mid \det(Q) = 1\}$\\
		$\Skew(n)$ & vector space of skew-symmetric matrices $\{A\in \R^{n\times n}\mid A^T = -A\}$.\\
		$\mathcal{M}$ & a Riemannian manifold\\
		$St(n,p)$ & Stiefel manifold $St(n,p)=\{U\in\R^{n\times p}\mid U^TU = I_p\}$\\
		$T_USt(n,p)$ & tangent space of $St(n,p)$ at $U\in St(n,p)$, 
		$T_USt(n,p)=\{\Delta \mid U^T\Delta +\Delta^TU =0\}$\\
		$\exp_m$ & matrix exponential $\exp_m(X) = \sum_{k=0}^\infty \frac{1}{k!}X^k$\\
		$\text{Exp}_U,\,\text{Log}_U$ & Riemannian exponential and Riemannian logarithm at $U\in St(n,p)$\\
		$i(\M)$ & injectivity radius of the Riemannian manifold $\M$\\
		$\mathcal{K}(\Delta_1,\Delta_2)$ & sectional curvature of the subspace spanned by the linear independent tangent vectors $\Delta_1,\Delta_2$\\
		$\text{Cut}(p)$ & set of cut points of all geodesics starting from $p\in\M$ 
	\end{tabular}
\end{small}

\section{Background theory} \label{chap:background}

    In this chapter, we recap basic concepts of manifolds and show how they apply to the Stiefel manifold. 
    \Cref{sec:background} introduces basic manifold theory, \Cref{sec:quotient} reviews the essentials of quotient spaces of Lie groups by Lie subgroups. In \Cref{sec:StiefelManifold}, the general concepts are substantiated for the Stiefel manifold.

	\subsection{Geometric concepts on manifolds}\label{sec:background}
	The following section follows the discussions from \cite{Zim21}, \cite{GeodesicSource} and \cite{leesmoothmanifold}. 
   A fundamental concept is that of a tangent space.

\begin{definition}[Tangent space]
\label{def:gentang}
Let $\M$ be a differentiable manifold. The tangent space of $\M$ at a point $p\in \M$ is defined as the space of velocity vectors of all differentiable curves $c\colon t\mapsto c(t)$ passing through $p$:
\begin{equation*}
T_p\M = \{\dot{c}(t_0)|c\colon J\to \M,\,c(t_0) = p\}.
\end{equation*}
Here $J\subseteq\R$ is an arbitrarily small open interval with $t_0\in J$.
\end{definition}
For embedded submanifolds $\M\subset \R^{n+d}$, the velocity vector $\dot c(t)\in\R^{n+d}$ is obtained by the usual rules of calculus. 
In the case of abstract manifolds, the symbol $\dot c(t)$ conceals quite a high level of abstraction. In absence of a surrounding space, the velocity vector of a curve is a derivative operator that induces a directional derivative of scalar functions. We omit the details.

The tangent bundle of $\M$ is the disjoint union of all tangent spaces
\[T\M = \{T_p\M|p\in\M\}.\]
Geometry begins with a scalar product for tangent vectors.
\begin{definition}[Length of a curve]
\label{def:lengthofcurve}
Let $\M$ be a differentiable manifold. 
A Riemannian metric on $\M$ is a family of inner products $\langle\cdot,\cdot\rangle_p\colon T_p\M \times T_p\M \to \R$, which is smooth in changes of the base point $p\in\M$.\\
The length of a tangent vector $v\in T_p \M$ is $\Vert v\Vert_p := \sqrt{\langle v,v\rangle_p}$. The length of a curve $c\colon\left[a,b\right]\to \M$ is defined as
\begin{equation*}
L(c) := \int_a^b \Vert\dot{c}(t)\Vert_{c(t)}dt.
\end{equation*}
The Riemannian distance between two points $p,\,q\in \M$ is 
\begin{equation*}
\dist_{\M}(p,q) := \inf\{L(c)\},
\end{equation*}
where $c$ is a piecewise smooth curve on the manifold connecting $p$ and $q$.
By convention, $\inf\{\emptyset\} = \infty$.
\end{definition}

Geodesics are candidates for length-minimizing curves and are characterized by the fact that they have no intrinsic acceleration.

\begin{definition}[Geodesics]
\label{def:geod2} 
A differentiable (unit speed) curve $\gamma\colon\left[0,1\right]\to \M$ is called geodesic (w.r.t. a given Riemannian metric) if the covariant derivative of the velocity vector field vanishes, i.e., 
\begin{equation} \label{eq:geodeq}
\frac{D\dot{\gamma}}{dt}(t) = D_t\dot{\gamma} = 0,\,\forall t\in\left[0,1\right]
\end{equation}
holds.
\end{definition}
Hence, geodesics are local solutions to an ordinary differential equation and depend smoothly on the initial values.
The Riemannian exponential map is based on these facts.

\begin{definition}[Riemannian exponential]
\label{def:riemanexp}
Let $\gamma_{p,v}$ be the geodesic starting from $p$ with velocity $v$. The Riemannian exponential is defined as
\begin{equation*}
\text{Exp}_p^{\M}\colon T_p \M\supset B_{\epsilon}(0)\to \M,\quad v\mapsto q := \gamma_{p,v}(1).
\end{equation*}
For technical reasons, $\epsilon>0$ must be small enough so that $\gamma_{p,v}(t)$ is defined on the unit interval $\left[0,1\right]$.
\end{definition} 

The Riemannian exponential is a local diffeomorphism \cite[Lemma 5.10]{GeodesicSource}.

\begin{definition}[Riemannian logarithm]
\label{def:riemanlog}
\label{def:Rielog}
The continuous inverse of the Riemannian exponential is called the Riemannian logarithm and is defined as
\begin{equation*}
\text{Log}_p^{\M}\colon \M\supset D_p\to B_{\epsilon}(0)\subset T_p \M,\quad q\mapsto v:=(\text{Exp}_p^{\M})^{-1}(q).
\end{equation*}
Here, $v$ satisfies $c_{p,v}(1) = q$.
\end{definition}

The size of the domain, where the Riemannian logarithm is well-defined is quantified by the injectivity radius.

\begin{definition}[Injectivity radius]\label{def:injrad}
Let $\epsilon$ be the maximum radius of $B_{\epsilon}(0)$ such that the Riemannian exponential at $p$, $\text{Exp}_p^{\M}\colon T_p \M\supset B_{\epsilon}(0)\to D_p\subset\M$, is invertible. 
Then, $\epsilon$ is called the injectivity radius of $\M$ at $p$ and is denoted by $i_p(\M)$.\\
The infimum of $i_p(M)$ over all $p\in\M$ is called injectivity radius of $\M$,
\[i(\M) = \inf_{p\in\M}i_p(\M).\]
\end{definition}
For a more detailed exploration of the notion of the injectivity radius, see \cite[Chap. 13]{docarmo}.\\
The  Riemannian exponential (tangent space to manifold) and the Riemannian logarithm (manifold to tangent space) 
form a special set of coordinate charts, called the Riemannian normal coordinates.
The normal coordinates are radially isometric in the sense that the Riemannian distance between $p$ and $q = \text{Exp}_p^{\M}(v)$ is the same as the length of the tangent vector $\Vert v\Vert_p = \Vert \text{Log}_p^{\M}(q)\Vert_p$.

\subsection{Quotient manifolds}\label{sec:quotient}
    Manifolds that arise as quotients of Lie groups by  Lie subgroups are highly structured.
    The Stiefel manifold belongs to this class. We recap the essentials of this quotient space construction and refer to \cite{gallier,leesmoothmanifold} for the details.
    
  A Lie group is a differentiable manifold $\mathcal{G}$ which also features a group structure, such that the group operations "multiplication" and "inversion" are both smooth. 
  Let $\mathcal{H}\le\mathcal{G}$ be a Lie subgroup and $p\in\mathcal{G}$. A subset of $\mathcal{G}$ of the form $\left[p\right] := p\mathcal{H} = \{p\cdot q|q\in\mathcal{H}\}$ is called left coset of $\mathcal{H}$. The left cosets form a partition of $\mathcal{G}$ and the quotient space determined by this partition is called the left coset space of $\mathcal{G}$ modulo $\mathcal{H}$ and is denoted by $\mathcal{G}/\mathcal{H}$. 

\begin{theorem}[cf. {\cite[Thm. 21.17]{leesmoothmanifold}}]\label{thm:quotman}
Let $\mathcal{G}$ be a Lie group and let $\mathcal{H}$ be a closed subgroup of $\mathcal{G}$. Then the left coset space $\mathcal{G}/\mathcal{H}$ is a manifold of dimension $\dim\mathcal{G} - \dim\mathcal{H}$ with a unique smooth structure such that the quotient map $\pi:\mathcal{G}\to\mathcal{G}/\mathcal{H}$, $p\mapsto\left[p\right]$ is a smooth submersion. The left action of $\mathcal{G}$ on $\mathcal{G}/\mathcal{H}$ given by 
\[p_1\cdot(p_2\mathcal{H}) = (p_1p_2)\mathcal{H}\]
turns $\mathcal{G}/\mathcal{H}$ into a homogeneous $\mathcal{G}$-space.
\end{theorem}

A homogeneous $\mathcal{G}$-space is a differentiable manifold endowed with a transitive smooth action by a Lie group $\mathcal{G}$. The fact that the action is transitive means that the structure "looks the same" everywhere on the manifold. 
Each preimage $\mathcal{G}_q:=\pi^{-1}(q)\subset \mathcal{G}$ is called fiber over $q$ and is itself a closed embedded submanifold. Let $\langle\cdot,\cdot\rangle_p^{\mathcal{G}}$ be the Riemannian metric of $\mathcal{G}$ at each point $p\in\mathcal{G}$. Then the tangent space $T_p\mathcal{G}$ decomposes into an orthogonal direct sum $T_p \mathcal{G} = T_p \mathcal{G}_{\pi(p)}\oplus (T_p \mathcal{G}_{\pi(p)})^{\perp}$ with respect to the metric. The tangent space of the fiber $V_p := T_p \mathcal{G}_{\pi(p)}$ is called the vertical space and is described by the kernel of the differential $d\pi_p\colon T_p \mathcal{G} \to T_{\pi(p)} \mathcal{G}/\mathcal{H}$. The orthogonal complement of the vertical space $V_p^{\perp} =: H_p$ is called the horizontal space. A crucial insight is that the tangent space of the quotient at $\pi(p)$ may be identified with the horizontal space at $p$, i.e.,
\begin{equation*}
    H_p\cong T_{\pi(p)} \mathcal{G}/\mathcal{H}.
\end{equation*}

\begin{remark}[{cf. \cite[p. 212]{ONeill1983}}]
\label{rem:quotisisom}
For every tangent vector $w\in T_{\pi(p)}\mathcal{G}/\mathcal{H}$ there is $\bar{x} = \bar{v} + \bar{w}\in V_p\oplus H_p = T_p\mathcal{G}$ such that $d\pi_p(\bar{x}) = w$. The horizontal component $\bar{w}\in H_p$ is unique and is called the horizontal lift of $w$. By relying on horizontal lifts, a Riemannian metric on the quotient can be defined by 
\begin{equation*}
\langle w_1,w_2\rangle_{\pi(p)}^{\mathcal{G}/\mathcal{H}} := \langle \bar{w_1},\bar{w_2}\rangle_p^\mathcal{G}
\end{equation*} 
for $w_1,w_2\in T_{\pi(p)}\mathcal{G}/\mathcal{H}$. With respect to this and only this metric, by construction, $d\pi|_{H_p}$ preserves the inner products of horizontal vectors and thus describes an isometry between the horizontal space $H_p$ and $T_{\pi(p)}\mathcal{G}/\mathcal{H}$. As a consequence, horizontal geodesics in $\mathcal{G}$ are mapped to geodesics on $\mathcal{G}/\mathcal{H}$ under $\pi$. Horizontal geodesics are geodesics in the total space whose velocity fields remain in the horizontal space for all time $t$.
\end{remark}

\subsection{The Stiefel Manifold}\label{sec:StiefelManifold}
In this section, we give a short introduction to the Stiefel manifold, see \cite{EdelmanAriasSmith1998, Zim21}.
Let $p\le n$: The set of all rectangular matrices with orthonormal columns 
\begin{equation*}
St(n,p) := \{U\in\R^{n\times p}|U^T U = I_p\}
\end{equation*}
is called (compact) Stiefel manifold.
The Stiefel manifold $St(n,p)$ is a submanifold of $\R^{np}$ of dimension 
\begin{equation}\label{eq:Stidim}
np - \frac{1}{2}(p(p+1)) = \frac{1}{2}p(p-1) + (n-p)p.
\end{equation}
The tangent space at a point $U\in St(n,p)$ is  
\begin{equation*}
T_{U}St(n,p) = \{\Delta\in\R^{n\times p}|U^T\Delta = -\Delta^T U\}.
\end{equation*}

Tangent vectors $\Delta$ can be represented in either of the following forms:
$\Delta = UA + (I - UU^T)T$, or $\Delta = UA + U^{\perp}B$,
where $A\in\Skew(p)$ and $T\in\R^{n\times p}$ and $B\in\R^{(n-p)\times p}$ arbitrary.

The Stiefel manifold is a quotient space of the orthogonal group $St(n,p) = O(n)/H (= O(n)/O(n-p))$, 
the associated left cosets are 
\begin{equation*}
\left[Q\right] = \left\{Q\begin{bmatrix} I_p & 0 \\ 0 & Q_{n-p}\end{bmatrix}\big|Q_{n-p}\in O(n-p)\right\}.
\end{equation*}
Therefore, the Stiefel manifold is a homogeneous $O(n)$-space, which implies that we can move from any point $U\cong\left[Q\right]\in St(n,p)$ to any other point $\tilde{U}\cong\left[\tilde{Q}\right]\in St(n,p)$ by left-multiplication with a certain $W\in O(n)$.\\
Let $Q\in O(n)$. The vertical space $V_Q$ and the horizontal space with respect to the (scaled) inner product $\langle X,Y \rangle := \frac12 \langle X,Y\rangle_F = \frac12 \trace(X^TY)$ on $O(n)$ are
\begin{equation*}
V_Q = \left\{Q\begin{bmatrix}0 & 0\\0 & C\end{bmatrix}\in\R^{n\times n}\big|\,C\in\Skew(n-p)\right\}.
\end{equation*}
and 
\begin{equation*}
H_Q = \left\{Q\begin{bmatrix}A & -B^T\\B & 0\end{bmatrix}\in\R^{n\times n}\big|\,A\in\Skew(p),\,B\in\R^{(n-p)\times p}\right\},
\end{equation*}
see \cite{EdelmanAriasSmith1998}.
Intuitively, motion in the direction of the vertical space brings no changes in quotient space $St(n,p)$, since the first $p$ columns are not touched. Therefore, concepts like metric and geodesic may be restricted to the horizontal space, which can be identified with the tangent space of the quotient manifold $St(n,p)$ as described in \Cref{sec:quotient}.\\
Matrices $\Delta_1,\,\Delta_2$ from the horizontal space $H_Q\cong T_{\left[Q\right]} St(n,p)$ 
can be considered both as tangent vectors of
the quotient space $St(n,p)$ and as special tangent vectors of the total space $O(n)$. We obtain an inner product for the former by recycling the inner product of the latter,
\begin{eqnarray*}
\langle\Delta_1,\Delta_2\rangle &=& \frac12 \trace\left(\begin{bmatrix} A_1^T & B_1^T\\ -B_1 & 0\end{bmatrix}Q^T Q\begin{bmatrix}A_2 & -B_2^T\\ B_2 & 0\end{bmatrix}\right)\\
&=& \frac12\trace(A_1^T A_2) + \trace(B_1^T B_2).
\end{eqnarray*}
This is called the canonical metric on $St(n,p)$, cf. {\cite[eq. 2.22]{EdelmanAriasSmith1998}}).
Under this metric, the geodesic that starts from $U\in St(n,p)$ with velocity $\Delta=UA + U^\bot B\in T_U St(n,p)$ reads 
\begin{equation}\label{eq:geodfull}
    \gamma(t) = \begin{bmatrix}U & U^{\perp}\end{bmatrix}\exp_m\left(t\begin{bmatrix}A&-B^T\\B&0\end{bmatrix}\right)\begin{bmatrix} I_p\\ 0\end{bmatrix}.
\end{equation}
For a detailed derivation see \cite{EdelmanAriasSmith1998}. The Riemannian exponential immediately emerges as $\text{Exp}_U(\Delta) = \gamma(1)$, wherever well-defined. 

For the Riemannian Logarithm, there does not exist a closed formula. The central objective is to find a tangent vector $\Delta\in T_U St(n,p)$ for two given points $U,\tilde U\in St(n,p)$ such that $\text{Exp}_U(\Delta) = \tilde U$. For simplicity, we speak of a tangent vector that connects $U$ and $\tilde U$. To give an idea on deriving the Riemannian Logarithm, we follow \cite{ralfcanonical}. Assume we have found a tangent vector $\Delta$ connecting $U,\tilde U$. Let $\Delta$ be parametrized by the matrices $A$ and $B$ and let $\exp_m\left(\begin{bsmallmatrix}A & -B^T \\ B & 0\end{bsmallmatrix}\right) = \begin{bsmallmatrix} M&X\\N&Y\end{bsmallmatrix}\in SO(n)$. By the Riemannian exponential, we obtain
\begin{align*}
    \tilde{U} = \text{Exp}_U(\Delta) = UM + U^\perp N.
\end{align*}
With only $U, \tilde U$ given, we immediately obtain $M = U^T\tilde U$ and $N = (U^\perp)^T\tilde U$. To recover the matrices $A$ and $B$ that determine the requested tangent vector $\Delta$, we need a suitable orthogonal completion $X,Y$ to $M,N$ such that $\log_m\left(\begin{bsmallmatrix}M&X\\N&Y\end{bsmallmatrix}\right) = \begin{bsmallmatrix}A & -B^T \\ B & 0\end{bsmallmatrix}$. It is confirmed by \cite[Thm. 3.1]{ralfcanonical} that $\Delta = UA + U^\perp B$ is a tangent vector connecting $U$ and $\tilde U$ if a suitable orthogonal completion is found, such that the matrix logarithm of $\begin{bsmallmatrix}M&X\\N&Y\end{bsmallmatrix}$ produces a zero in the lower right block. Hence, the task of finding a tangent vector connecting two points on the manifold boils down to finding a rotation $\Phi\in SO(p)$ such that $X_0\Phi,Y_0\Phi$ is a suitable orthogonal completion to $M,N$. Here, $X_0,Y_0$ is some arbitrary orthogonal completion. An algorithm for the Riemannian Logarithm in given by \cite[Algorithm 3.1]{ralfcanonical}.

	\section{On Cut Points and a bound of the Injectivity Radius}\label{sect:invoninj}
	
	In this section, we re-derive Rentmeesters' bound \cite{Quentin} on the injectivity radius of the Stiefel manifold,
	\[i(St(n,p)) \ge \sqrt{\frac{4}{5}}\pi\approx 2.81.\]
	It is an open question, whether this bound is sharp.
	\footnote{Note. On March 4, 2024, two days before the submission of the first preprint of this work, personal communication revealed that Absil/Mataigne worked independently on the injectivity radius of the Stiefel manifold, but for a parametric family of Riemannian metrics \cite{HueperMarkinaLeite2020, log_comp, nguyen2022curvature}.
		By the discovery of conjugate points and geodesic loops on the Stiefel manifold, an upper bound on the injectivity radius is obtained. Their work is now available as a preprint \cite{absil2024ultimate}. } 
	\\
	In the following, we recap the theoretical foundation for working with cut points and conjugate point in relation to the injectivity radius. Our main references are \cite{docarmo,petersen2016riemannian,gallier} and \cite{Quentin}. \\
	We start by defining the sectional curvature of a tangent plane section of a manifold $\M$. 
	Given a vector space $V$, we write
	\[
	\Vert X\wedge Y\Vert:=\sqrt{\Vert X\Vert^2\Vert Y\Vert^2 - \langle X,Y\rangle^2}, \text{ for } X,Y\in V.
	\]

	\begin{definition}[cf. {\cite[Chapter 4, Prop. 3.1, Def. 3.1]{docarmo}}]
		Let $\Omega\subset T_p\M$ be a two-dimensional subspace of the tangent space $T_p\M$ and let $X,Y\in\Omega$ be two linearly independent vectors. Then
		\[\mathcal{K}_p(\Omega) := \mathcal{K}_p(X,Y) := \frac{\langle R(X,Y)X,Y\rangle}{\Vert X\wedge Y\Vert^2}\]
		is called the sectional curvature of $\Omega$ at $p$. 
		It does not depend on the choice of the basis vectors $X,Y\in\Omega$.
	\end{definition}
	
	Here $R$ denotes the curvature tensor of $\M$. For details see \cite{docarmo}. 
	An explicit formula for determining the sectional curvature of the Stiefel manifold is given in \cite[Prop. 4.2, eq. (34)]{nguyen2022curvature}. Next, we introduce the notion of Jacobi fields. 
	
	\begin{definition}[Jacobi Field (cf. {\cite[Chapter 5, Def. 2.1]{docarmo}})] 
		Let $\gamma\colon\left[0,1\right]\to\M$ be a geodesic in $\M$. A vector field $J$ along $\gamma$, i.e., $J(t)\in T_{\gamma(t)}\M$, is said to be a Jacobi field if it satisfies the Jacobi equation
		\begin{equation}\label{eq:jacobieq}
			\frac{D^2 J}{dt^2} + R(\gamma^\prime(t),\,J(t))\gamma^\prime(t) = 0,
		\end{equation}
		for all $t\in\left[0,1\right]$.
	\end{definition}
	
	A Jacobi field is determined by the initial conditions $J(0)$ and $\frac{DJ}{dt}(0)$. If the dimension of $\M$ is $m$, there are $2m$ linearly independent Jacobi fields along a geodesic $\gamma$. 
	Specifying the initial condition $J(0) = 0$, we are left with $m$ linearly independent Jacobi fields along $\gamma$. 
	The next lemma gives a characterization for Jacobi fields with $J(0) = 0$.
	
	\begin{lemma}[cf. {\cite[Prop 17.22]{gallier}}]\label{col:jacobialt}
		Let $\gamma\colon\left[0,1\right]\to\M$ be a geodesic and let $W\in T_{\gamma^\prime(0)}(T_{\gamma(0)}\M)\cong T_{\gamma(0)}\M$. Then a Jacobi field $J$ along $\gamma$ with $J(0) = 0$ and $\frac{DJ}{dt}(0) = W$ is given by
		\begin{align*}
			J(t) = (d\text{Exp}_{\gamma(0)})_{t\gamma^\prime(0)}[tW],\quad t\in\left[0,1\right].
		\end{align*}
	\end{lemma}
	
	With the Jacobi fields at hand, conjugate points can be defined.
	
	\begin{definition}[Conjugate Point (cf. {\cite[Chapter 5, Def. 3.1]{docarmo}})]
		Let $\gamma\colon\left[0,1\right]\to\M$ be a geodesic. The point $\gamma(t_0)$ is said to be conjugate to $\gamma(0)$ along $\gamma$, if there exists a non trivial Jacobi field $J$ along $\gamma$ with $J(0) = 0 = J(t_0)$.\\
		The maximum number of such linearly independent fields is called the multiplicity of the conjugate point $\gamma(t_0)$.
	\end{definition}
	
	A conjugate point to $\gamma(0)$ can be identified with a critical point of the Riemannian exponential $\text{Exp}_{\gamma(0)}$, see \cite{docarmo}. 
	Classical Riemannian geometry provides a statement about the distance between conjugate points.
	
	\begin{proposition}[cf. {\cite[Theorem 6.4.6]{petersen2016riemannian}}]\label{prop:distconj}
		Let  $\M$ be a Riemannian manifold. Suppose that for any $p\in\M$, the sectional curvatures are bounded by $\mathcal{K}_p\le H$ with $H>0$ being a constant.
		Then
		\[
		\mathrm{Exp}_p:B(0,\pi/\sqrt{H})\to\M
		\]
		has no critical, hence no conjugate points.
	\end{proposition}
	
	The last term we introduce, before a bound of the injectivity radius can be formulated, is that of a cut point of a geodesic. 
	Let $\M$ be a complete Riemannian manifold in the following (a property featured by the Stiefel manifold).
	
	\begin{definition}[Cut Point (cf. {\cite[p. 267]{docarmo}})]
		Let $\M$ be a complete Riemannian manifold, let $p\in\M$ and let $\gamma\colon\left[0,\infty\right)\to\M$ be a normalized geodesic with $\gamma(0) = p$. We know that if $t>0$ is sufficiently small, $d(\gamma(0),\gamma(t)) = t$, i.e., $\gamma\left(\left[0,t\right]\right)$ is a minimizing geodesic (see \cite[Chapter 3, Prop. 3.6]{docarmo}). In addition, if $\gamma\left(\left[0,t_1\right]\right)$ is not minimizing, the same is true for all $t>t_1$ (see \cite[Prop. 16.18]{gallier}). By continuity, the set of numbers $t>0$ for which $d(\gamma(0),\gamma(t)) = t$ is of the form $\left[0,t_0\right]$ or $\left[0,\infty\right)$. In the first case, $\gamma(t_0)$ is called the cut point of $p$ along $\gamma$. In the second case, we say that such cut point does not exist.\\
		We define the cut locus of $p$, denoted by $\text{Cut}(p)$, as the union of the cut points of $p$ along all geodesics starting from $p$. 
	\end{definition}
	
	So, a cut point of a geodesic can be seen as the location from which on the geodesic fails to describe a unique shortest path. A fundamental property of cut points is the following.
	
	\begin{proposition}[cf. {\cite[Chapter 13, Prop. 2.2]{docarmo}}] \label{prop:cutpoint}
		Suppose that $\gamma(t_0)$ is the cut point of $p = \gamma(0)$ along $\gamma$. Then
		\begin{enumerate}
			\item either $\gamma(t_0)$ is the first conjugate point of $\gamma(0)$ along $\gamma$,
			\item or there exists a geodesic $\sigma\neq\gamma$ from $p$ to $\gamma(t_0)$ such that $L(\sigma)=L(\gamma)$.
		\end{enumerate}
		Conversely, if one of the above conditions is met, then there exists $\tilde{t}$ in $\left(0,t_0\right]$ such that $\gamma(\tilde{t})$ is the cut point of $p$ along $\gamma$.
	\end{proposition}
	
	\begin{corollary}[cf. {\cite[Chapter 13, Cor. 2.8]{docarmo}}]\label{cor:qsmin}
		If $q\in\M\setminus \text{Cut}(p)$, there exists a unique minimizing geodesic joining $p$ to $q$. 
	\end{corollary} 
	
	Another version of this corollary can be found in \cite[Thm. 17.30]{gallier}.\\ 
	\Cref{cor:qsmin} shows that $\text{Exp}_p$ is injective on an open ball $B_r(p)$ if and only if the radius $r$ is less than or equal to the distance from $p$ to $\text{Cut}(p)$. For this reason, we can write 
	\[i(\M) = \inf_{p\in\M}d(p,\text{Cut}(p))\]
	for the injectivity radius of $\M$, see \cite[p. 271]{docarmo}.
	
	\begin{proposition}[cf. {\cite[Prop. 17.32b]{gallier}}]\label{prop: closedgeod}
		For $p\in\M$ suppose $q\in \text{Cut}(p)$ realizes the distance from $p$ to $\text{Cut}(p)$, i.e., $d(p,q)= d(p,\text{Cut}(p)) =: l$. If there are no minimal geodesics from $p$ to $q$ such that $q$ is conjugate to $p$ along this geodesic, then there are exactly two minimizing geodesics $\gamma$ and $\sigma$ from $p$ to $q$, with $\sigma^\prime(l) = -\gamma^\prime(l)$ and $l = d(p,q)$. If, in addition, $d(p,q) = i(\M)$, then $\gamma$ and $\sigma$ together form a closed geodesic. 
	\end{proposition} 
	
	With the above preparations, we are now in a position to state the classical bound for the injectivity radius. 
	
	\begin{theorem}[Klingenberg, stated as Lemma 6.4.7 in \cite{petersen2016riemannian}]
		\label{thm:Klingenberg}
		Let $(\M,\langle\cdot,\cdot\rangle)$ be a compact Riemannian manifold with sectional curvatures bounded by $\mathcal{K}_p\leq H$, where $H>0$. Then
		the injectivity radius at any $p\in \M$ satisfies
		\[
		i(p)\geq \min\left\{\frac{\pi}{\sqrt{H}}, \frac12 l_p\right\},
		\]
		where $l_p$ is the length of a shortest closed geodesic starting from $p$.
		For the global injectivity radius, it holds
		\[
		i(\M)  \geq\frac{\pi}{\sqrt{H}}
		\quad \text{ or } \quad \mathrm{inj}(\M) =\frac12 l,
		\]
		where $l$ is the length of a shortest closed geodesic on $\M$.
	\end{theorem}

	Notice the correspondence between this statement and items 1. and 2. of \Cref{prop:cutpoint}.\\
	
	Applying the theoretical framework to the Stiefel manifold, one obtains a concrete bound on its injectivity radius. 
	In \cite{ZimmermannStoye2024},  a global bound on the sectional curvature is given
	\[\mathcal{K}_p\le \frac{5}{4} =: H.\]
	In particular, this bound is sharp for $p\geq 2, n\geq p+2$ and is only achieved for the tangent plane spanned by the normalized, orthogonal tangent vectors $\Delta_1 = UA_1 + U^{\perp}B_1$, $\Delta_2 = UA_2 + U^{\perp}B_2$ with\\
	\[A_1 = A_2 = 0,\]
	\[B_1 = \frac{1}{\sqrt{2}}\left[\begin{array}{c|c}
		\begin{matrix} 0 & 1\\ 1 & 0 \end{matrix} & \mathbf{0} \\
		\hline
		\mathbf{0}& \mathbf{0} \\ 
	\end{array}\right]
	\mbox{ and }
	B_2 = 
	\frac{1}{\sqrt{2}}\left[\begin{array}{c|c}
		\begin{matrix} 1  & 0\\ 0 & -1 \end{matrix} & \mathbf{0} \\
		\hline
		\mathbf{0}& \mathbf{0} \\ 
	\end{array}\right].\]
	Here, the matrices $B_1$ and $B_2$ are unique up to trace-preserving transformations (see\cite{ZimmermannStoye2024}).\\
	An explicit formula for the calculation of the sectional curvature on the Stiefel manifold is given in \cite[Prop. 4.2, eq. (34)]{nguyen2022curvature}. The calculation only depends on the parametrization of the tangent vectors via $A$ and $B$ and not on the point $U$. Hence, we write $\mathcal{K}$ instead of $\mathcal{K}_p$ in the following.\\
	
	With \Cref{thm:Klingenberg}, it now either holds 
	\begin{equation}\label{eq:injrad}
		i(St(n,p))\ge\sqrt{\frac{4}{5}}\pi\approx 2.81,
	\end{equation}
	or there exists a closed geodesic, with length less than $2\sqrt{\frac{4}{5}}\pi$. In \cite[p. 94]{Quentin} and also in \cite[Theorem 6.1]{absil2024ultimate}, however, it is shown that closed geodesics in $St(n,p)$ have at least length $2\pi$. Thus, for the injectivity radius of the Stiefel manifold, the estimate \eqref{eq:injrad} holds.\\
	A special feature of the Stiefel manifold as a homogeneous $O(n)$-space is described in the following. Let $\gamma:\left[0,1\right]\to St(n,p)$ be a geodesic starting from $U$ with $\Delta = UA + U^\perp B$. Further, let $Q\in O(n)$ be an arbitrary orthogonal matrix and define $\hat \gamma := Q\gamma$. Then $\hat\gamma$ defines the geodesic with starting point $Q U$ and velocity $Q^T\Delta = (QU)A + (QU^\perp)B$. Hence, both tangents are parameterized by the same matrices $A$ and $B$. Since the calculation of the length of the geodesics and the evaluation of Jacobi fields as the directional derivative of the Riemannian exponential depends only on the matrices $A$ and $B$ and not on the starting points, both geodesics have the same cut point. 
	Therefore, it holds $d(U,\text{Cut}(U)) = d(\tilde U,\text{Cut}(\tilde U))$ for any two points $U,\tilde U\in St(n,p)$ of the Stiefel manifold. As a result, the following corollary emerges.
	
	\begin{corollary}\label{cor:InjRadOnePoint}
		Let $U\in St(n,p)$, then it holds 
		\[i(St(n,p)) = d(U,\text{Cut}(U)).\]
	\end{corollary}
	
	Furthermore, we can use the second case in the \Cref{prop:cutpoint} to give a limitation for cut points, which are no conjugate points along any geodesic. 
	
	\begin{remark}\label{rem:2piinsecond}
		If we exclude all conjugate points from the consideration of cut points, this new "injectivity radius", $\hat{i}$, is described by the shortest closed geodesic $\hat{\gamma}$, i.e., $\hat{i} = \frac{1}{2}L(\hat{\gamma})$ (see \Cref{thm:Klingenberg}). This exclusion of conjugate points is possible because the proof of \Cref{prop:cutpoint} clearly distinguishes between conjugate points and no conjugate points.\\
		On the Stiefel manifold the shortest closed geodesics have at least length $2\pi$. In particular, the shortest closed geodesics have exactly the length $2\pi$. An example of a closed geodesic of length $2\pi$ is given by $A = 0$ and $B = \begin{bmatrix}2\pi&0\end{bmatrix}$.
		Thus, it follows with \Cref{prop:cutpoint} that cut points of geodesics on the Stiefel manifold which are not conjugate points of the geodesic are at least at a distance of $\frac{1}{2}2\pi = \pi$ along the geodesic. 
		Hence, the search for closest conjugate points is essential when investigating the bound \eqref{eq:injrad}.
	\end{remark}

		\section{Numerical experiments on sharpness of the injectivity radius bound}\label{sec:numexpInj}
\begin{figure}[h!]
  \centering
  \includegraphics[width=0.5\textwidth]{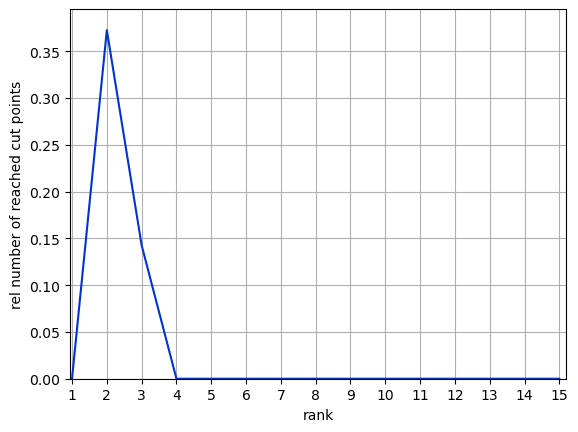}
  \caption{Number of reached cut points relative to the number of examined geodesics for different ranks of tangent vectors $\Delta$. 
  (variables: $n\in\left[4,100\right]$, $p\in\left[2,15\right]$ and $\mu\in\left[2.8,3.1\right]$).}
  \label{fig:ranks}
\end{figure}

\begin{figure}[h!]
  \centering
  \includegraphics[width=0.75\textwidth]{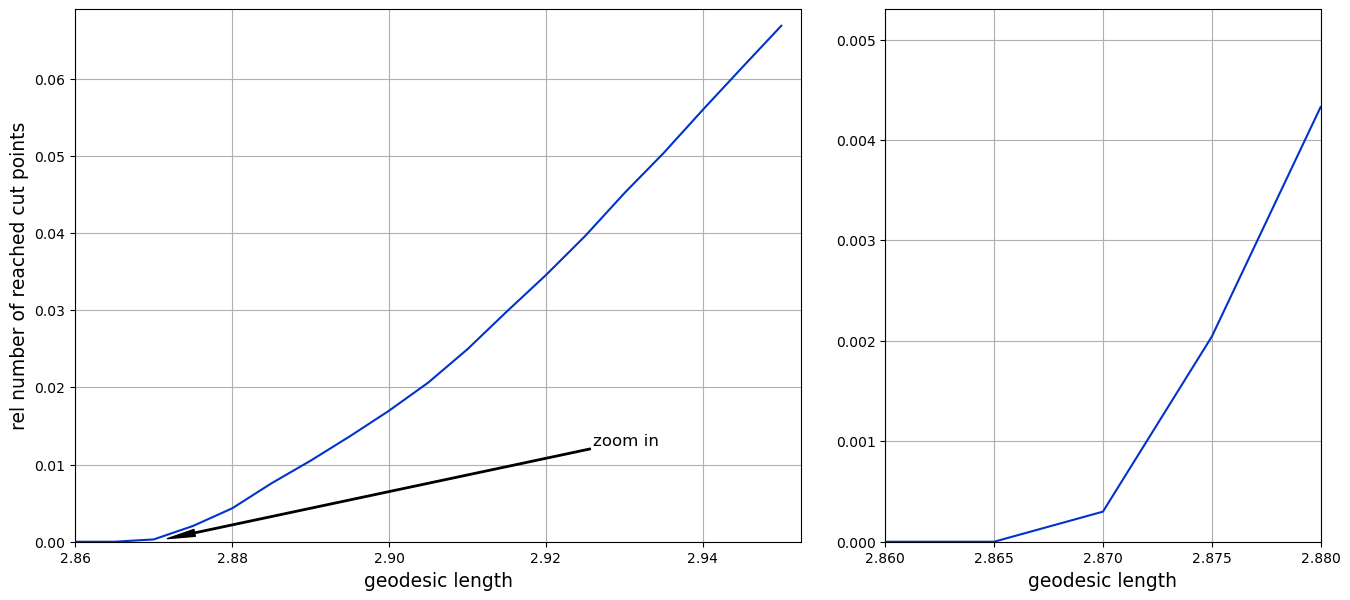}
  \caption{Number of reached cut points relative to the number of examined geodesics for different lengths $\mu$ and a zoom in on the right.\\ 
  (variables: $n\in\left[4,100\right]$, $p\in\left[2,7\right]$ and $\mu\in\left[2.86,2.95\right]$ discretized in step sizes of $0.005$).}
  \label{fig:mugraph}
\end{figure}

\begin{figure}[h!]
  \centering
  \includegraphics[width=0.5\textwidth]{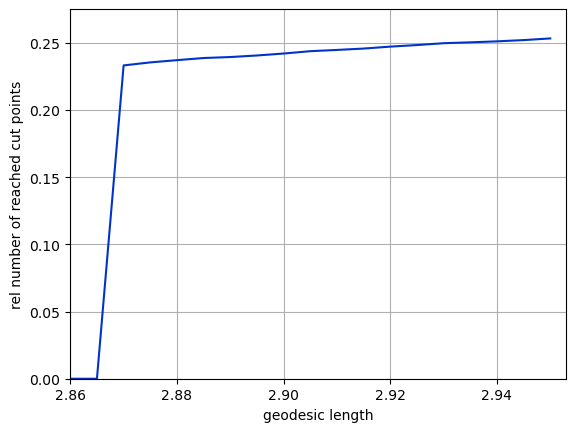}
  \caption{Number of reached cut points relative to the number of examined geodesics for different lengths $\mu$. The geodesics are defined by starting velocities from maximal sectional curvature planes.\\
  (variables: $n\in\left[4,100\right]$, $p\in\left[2,7\right]$ and $\mu\in\left[2.86,2.95\right]$ discretized in step sizes of $0.005$).}
  \label{fig:mugraph_maxcurv}
\end{figure}

The injectivity radius of the Stiefel manifold is equal to the distance of any point $U\in St(n,p)$ to the set of cut points of all geodesics starting from $U$ (see \Cref{cor:InjRadOnePoint}). We investigate the injectivity radius by examining geodesics of various lengths $\mu$ for their cut points. The geodesics are defined by the starting point $\begin{bsmallmatrix}I_p\\0\end{bsmallmatrix}$ and a `random' tangent vector $\Delta$. For each geodesic, we solve the geodesic endpoint problem for the endpoints of the geodesic to find a (possibly) different geodesic connecting the same points. If the geodesic found by solving the geodesic endpoint problem is shorter than the examined geodesic, then the examined geodesic is exposed as non-minimizing and therefore its cut point was already reached. This results in the injectivity radius being smaller than the length of the examined geodesic. In this way we approach the injectivity radius from above.\\
In first experiments with the geodesic length $\mu$ in the range $\left[2.8,\,3.1\right]$, it is noticed that especially for tangents of low rank, especially rank two, the cut points of the corresponding geodesics are already reached in the given interval (see \Cref{fig:ranks}). 
This may partly be explained by the fact that the maximum of the sectional curvature is only reached for tangent vectors of rank two and the limited range of directions in low dimensions. With larger sectional curvature the boundary for cut points, respectively conjugate points, is smaller, which creates the possibility for cut points of geodesics of smaller length.\\
Based on this, we restrict ourselves in the following investigations to tangents of low rank by choosing a small value for $p$. This is in line with the findings in \cite{ZimmermannStoye2024}.\\
In all further investigations with $\mu\in\left[2.8,\,3.1\right]$, the smallest geodesic length $\mu$ prior to which a geodesic had its cut point was found to be $2.87$. The interval is discretized in step sizes of $0.01$. To determine a more detailed statement about a bound for the injectivity radius, the shorter interval $\left[2.86,\,2.95\right]$ discretized with a step size of $0.005$ is examined (see \Cref{fig:mugraph}). Again, the smallest $\mu$ prior to which a geodesic had its cut point is $2.87$.\\
With a different approach, we further restrict the space of the analysed geodesics. From \Cref{thm:Klingenberg}, we know that the bound on the injectivity radius depends on a global constant $H$ that bounds the sectional curvature on the Stiefel manifold. The larger the sectional curvature, the smaller the bound on the injectivity radius and the smaller the lower bound on the length of geodesics possibly having cut points, respectively. Therefore, we now restrict ourselves to the study of `random' geodesics with starting velocity $\Delta$ coming from a tangent plane section of maximum sectional curvature $\frac54$. Tangent vectors spanning tangent plane sections with maximum sectional curvature can be found in \cite{ZimmermannStoye2024}. Following this approach, we obtain the same results as in the previous experiments (see \Cref{fig:mugraph_maxcurv}).\\
For Stiefel manifolds $St(n,2)$ with $p=2$ there is another method to investigate if there are other geodesics, and, in particular, shorter geodesics to certain points. Recall \Cref{sec:StiefelManifold}, where it is stated that finding a geodesic connecting two points boils down to finding a suitable rotation $\Phi\in SO(p)$. Let $U_0,\,U$ be points on the manifold (close enough to each other) and let $M=U_0^T U$ and $N = (U_0^\perp)^TU$. Moreover, let $\begin{bsmallmatrix}X_0\\Y_0\end{bsmallmatrix}$ be an orthogonal completion of $\begin{bsmallmatrix}M\\N\end{bsmallmatrix}$ such that $\det\left(\begin{bsmallmatrix}M&X_0\\N&Y_0\end{bsmallmatrix}\right) = 1$. Now let $\Phi\in SO(2)$ be such that
\[\text{log}_m\Big(\underbrace{\begin{bmatrix} M&X_0\Phi\\N&Y_0\Phi \end{bmatrix}}_{=:V_{\Phi}}\Big) = \begin{bmatrix} A&-B^T\\B&0 \end{bmatrix}\]
holds. Then $\text{Exp}_{U_0}(\Delta) = U$ with $\Delta = U_0A+QB\in T_{U_0} St(n,2)$. The matrices $\Phi\in SO(2)$ can be represented explicitly by
\[\Phi(\alpha) = \begin{bmatrix}\cos(\alpha)&-\sin(\alpha)\\\sin(\alpha)&\cos(\alpha) \end{bmatrix},\qquad\text{for $\alpha\in\left[-\pi,\,\pi\right)$}.\]
Therefore, the matrix solving the geodesic endpoint problem only depends on the scalar parameter $\alpha$. We are now able to iterate over this parameter to determine $\Phi(\alpha)$ such that 
\[\text{log}_m(V_{\Phi(\alpha)}) = \begin{bmatrix}A&-B^T\\B&0\end{bmatrix}.\]
In this way it is possible to determine different geodesics connecting the same points (if several exist). The procedure to investigate the injectivity radius is as follows. We start with geodesics of length $3.2$. To the endpoints of this geodesic, further geodesics connecting these points are determined using the procedure described above. As soon as a shorter geodesic is found, the length of the next investigated geodesics is reduced by 0.01. \\
\begin{figure}[h!]
  \centering
  \includegraphics[width=0.5\textwidth]{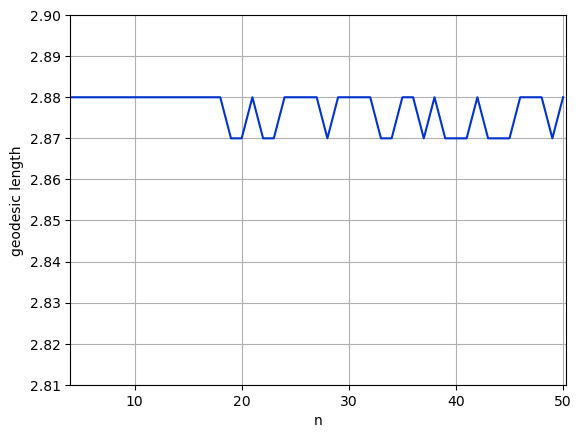}
  \caption{Smallest length $\mu$ of a geodesic to which endpoints a shorter geodesic is found\\ (variables: $n\in\left[4,50\right]$, $p = 2$, $\mu$ started at $3.2$, stopped after $250$ iterations)}
  \label{fig:stn2mus}
\end{figure}
In these studies, too, the smallest length $\mu$ of geodesics to whose endpoints there are shorter geodesics is at $2.87$ (see \Cref{fig:stn2mus}).\\
Summarizing the numerical experiments, it was not possible to reach the curvature-based bound of $\sqrt{4/5}\pi$ for the injectivity radius by means of the investigation of random geodesics. Rather, the experiments suggest that it is the conjugate points that determine the injectivity radius and not the curvature bound, and that the true value is around $2.87$. This coincides with the results for the canonical metric from the preprint \cite{absil2024ultimate} that appeared on arXiv two days before the first preprint of this work was submitted.  In case of the canonical metric, the conjecture made by the authors of \cite{absil2024ultimate} is that the injectivity radius is at $\sqrt{2}t_0$, where $t_0$ is the smallest positive root of $\frac{\sin(t)}{t} + \cos(t)$. Up to ten digits, $\sqrt{2}t_0 = 2.8690968494$.
In the next \Cref{sect:explcutpoint}, we construct an explicit example of a cut point of a geodesic at exactly this length and we find the same defining equation.
While the authors of \cite{absil2024ultimate} looked explicitly for critical points of the Riemannian exponential, we look for conjugate points. This is of course only a conceptual difference, since critical points and conjugate points are equivalent.

		\section{Explicit construction of a cut point under the canonical metric}\label{sect:explcutpoint}
In this section, we restrict ourselves to geodesics on the Stiefel manifold $St(4,2)$ starting in $U_0 = \begin{bsmallmatrix}I_2\\ 0\end{bsmallmatrix}$. 
This is the Stiefel manifolds of smallest dimension that features the maximum sectional curvature.
Obviously, $4\times 2$-Stiefel matrices are embedded in Stiefel manifolds of larger dimensions by just filling them up with zeros in a suitable way.
By \Cref{prop:cutpoint}, we know that cut points of geodesics are either first conjugate points or points to which there exist two different geodesics of equal length. Furthermore, by \Cref{rem:2piinsecond}, we know that the length of geodesics that feature a cut point which is not a conjugate point along the geodesic is at least $\pi$. In order to determine conjugate points along a geodesic $\gamma$, we investigate Jacobi fields 
\[J(t) = (d\text{Exp}_{U_0})_{t\gamma^\prime(0)}(tW)\]
along the geodesic for arbitrary directions $W\in T_{\gamma^\prime(0)}(T_{U_0}St(n,p))$. 
This requires us to calculate the directional derivative of the Stiefel exponential of $U_0$ at $t\gamma^\prime(0)$ in a direction $tW$. An explicit way to do this on Stiefel is outlined in \cite[Section 4.2]{ZimmermannHermite_2020}.\\
Let the tangent $\gamma^\prime(0)$ be parameterized by $A$ and $B$. The tangent space to a vector space can be identified with the vector space itself (see e.g. \cite[p.13]{grassmann}),  $T_{\gamma^\prime(0)}(T_{U_0}St(n,p))\cong T_{U_0}St(n,p)$. Therefore, let a direction $W\in T_{\gamma^\prime(0)}(T_{U_0}St(n,p))$ be parameterized analogous by $A_w$ and $B_w$. 

According to the definition of the Stiefel exponential, the calculation of its directional derivative boils down to extracting the first two columns from the directional derivative of the matrix exponential at $t\begin{bsmallmatrix} A & -B^T\\B&0\end{bsmallmatrix}$ in the direction $t\begin{bsmallmatrix} A_w & -B_w^T\\B_w&0\end{bsmallmatrix}$. Najfeld and Havel \cite{NAJFELD1995321} provide a formula for calculating the directional derivative of the matrix exponential. The directional derivative $D\exp_m(X)[Y]$ can be calculated via
\begin{align}\label{eq:expmDirectDeriv}\exp_m\begin{bmatrix}X&Y\\0&X\end{bmatrix} = \begin{bmatrix}\exp_m(X)&D\exp_m(X)[Y]\\0&\exp_m(X)\end{bmatrix}.\end{align}
This goes by the name of Mathias' Theorem in \cite[Theorem 3.6]{Higham:2008:FM}.
The geodesic $\gamma$ that we are going to investigate is defined by $U_0$ and a tangent vector $\Delta$ from a tangent plane section with maximum sectional curvature. We define
\[\Delta = U_0\underbrace{\begin{bmatrix}0&0\\0&0\end{bmatrix}}_{=:A} + U_0^{\perp}\underbrace{\begin{bmatrix} \frac{1}{2}&\frac{1}{2}\\-\frac{1}{2}&\frac{1}{2} \end{bmatrix}}_{=:B}.\]
By a standard result on Jacobi fields, see e.~g. \cite{docarmo}, there are $5 = \text{dim}(St(4,2))$ linearly independent Jacobi fields along $\gamma$ with $J(0) = 0$. Those can be found by calculating the Jacobi fields for linearly independent $W$, see \cite[Chapter 5, Remark 3.2]{docarmo}. Five linearly independent directions $W$ defining the Jacobi fields are 
\begin{align*}
    W_1 &= U_0\begin{bmatrix}0&-1\\1&0\end{bmatrix},\quad W_2 = U_0^\perp\begin{bmatrix}1&-1\\1&1\end{bmatrix}, \quad  W_3 = U_0^\perp\begin{bmatrix}1&1\\-1&1\end{bmatrix},\\
    W_4 &= U_0^\perp\begin{bmatrix}1&1\\1&-1\end{bmatrix},\quad 
    W_5 = U_0^\perp\begin{bmatrix}-1&1\\1&1\end{bmatrix}. 
\end{align*}
Calculating the Jacobi fields via the formula for the directional derivative of the matrix exponential \eqref{eq:expmDirectDeriv} can be done with the Jordan Canonical form. We obtain 
\begin{align*}
    J_1(t) &=(d\text{Exp}_{U_0})_{t\gamma^\prime(0)}[tW_1]= \begin{bmatrix}
        0&-\frac{t\cos\left(\frac{t}{\sqrt{2}}\right) + \sqrt{2}\sin\left(\frac{t}{\sqrt{2}}\right)}{2}\\
        \frac{t\cos\left(\frac{t}{\sqrt{2}}\right) + \sqrt{2}\sin\left(\frac{t}{\sqrt{2}}\right)}{2} & 0\\
        \frac{t}{2\sqrt{2}}\sin\left(\frac{t}{\sqrt{2}}\right)&-\frac{t}{2\sqrt{2}}\sin\left(\frac{t}{\sqrt{2}}\right)\\
        \frac{t}{2\sqrt{2}}\sin\left(\frac{t}{\sqrt{2}}\right)&\frac{t}{2\sqrt{2}}\sin\left(\frac{t}{\sqrt{2}}\right)
    \end{bmatrix},\\
    J_2(t) &=(d\text{Exp}_{U_0})_{t\gamma^\prime(0)}[tW_2]= \begin{bmatrix}
        0&0\\
        0&0\\
        \sqrt{2}\sin\left(\frac{t}{\sqrt{2}}\right)&-\sqrt{2}\sin\left(\frac{t}{\sqrt{2}}\right)\\
        \sqrt{2}\sin\left(\frac{t}{\sqrt{2}}\right)&\sqrt{2}\sin\left(\frac{t}{\sqrt{2}}\right)
        \end{bmatrix},\\
    J_3(t) &=(d\text{Exp}_{U_0})_{t\gamma^\prime(0)}[tW_3]= \begin{bmatrix}
        -\sqrt{2}t\sin\left(\frac{t}{\sqrt{2}}\right)&0\\
        0&-\sqrt{2}t\sin\left(\frac{t}{\sqrt{2}}\right)\\
        t\cos\left(\frac{t}{\sqrt{2}}\right)&t\cos\left(\frac{t}{\sqrt{2}}\right)\\
        -t\cos\left(\frac{t}{\sqrt{2}}\right)&t\cos\left(\frac{t}{\sqrt{2}}\right)
    \end{bmatrix},\\
    J_4(t) &=(d\text{Exp}_{U_0})_{t\gamma^\prime(0)}[tW_4]= \begin{bmatrix}
        0&-\sqrt{2}t\sin\left(\frac{t}{\sqrt{2}}\right)\\
        -\sqrt{2}t\sin\left(\frac{t}{\sqrt{2}}\right)&0\\
        t\cos\left(\frac{t}{\sqrt{2}}\right)&t\cos\left(\frac{t}{\sqrt{2}}\right)\\
        t\cos\left(\frac{t}{\sqrt{2}}\right)&-t\cos\left(\frac{t}{\sqrt{2}}\right)
    \end{bmatrix},\\
    J_5(t) &=(d\text{Exp}_{U_0})_{t\gamma^\prime(0)}[tW_5]= \begin{bmatrix}
        \sqrt{2}t\sin\left(\frac{t}{\sqrt{2}}\right)&0\\
        0&-\sqrt{2}t\sin\left(\frac{t}{\sqrt{2}}\right)\\
        -t\cos\left(\frac{t}{\sqrt{2}}\right)&t\cos\left(\frac{t}{\sqrt{2}}\right)\\
        t\cos\left(\frac{t}{\sqrt{2}}\right)&t\cos\left(\frac{t}{\sqrt{2}}\right)
    \end{bmatrix}.
\end{align*}
Those Jacobi fields are linearly independent. Moreover, the matrix function values $J_k(t)$ at any $t$ are linearly independent as long as all $t$-dependent entries are non-zero.\\
There is a conjugate point that can be extracted directly from the five Jacobi fields. It is obvious that $J_2$ vanishes for $t = \sqrt{2}\pi$. Therefore, $\gamma$ has a conjugate point at $\sqrt{2}\pi$. Next, we investigate whether there are linear combinations of Jacobi fields that define conjugate points that are closer to $U_0$. Because the matrices defined by the Jacobi fields at any $t$ are linearly independent as long as all $t$-dependent entries are non-zero, we are looking for linear combinations of Jacobi fields that vanish at some $t$ where a $t$-dependent entry becomes zero. The term $\sin\left(\frac{t}{\sqrt{2}}\right)$ becomes zero for multiples of $\sqrt{2}\pi$. The term $\cos\left(\frac{t}{\sqrt{2}}\right)$ becomes zero for $\sqrt{2}\left(\frac{\pi}{2} + n\pi\right)$, $n\in\Z$. The smallest positive root of the term $\frac12\left(t\cos\left(\frac{t}{\sqrt{2}}\right) + \sqrt{2}\sin\left(\frac{t}{\sqrt{2}}\right)\right)$ is $t_1$, where $t_1\approx 2.8690968494$. As the length of geodesics that feature conjugate points is bounded by the injectivity radius, the smallest candidate fulfilling $t\geq\sqrt{\frac{4}{5}}\pi$ to look for a linear combination of Jacobi fields vanishing at $t$ is $t_1$. In fact, $J(t) = J_1(t) - \frac{t_1}{4}J_2(t)$ is a Jacobi field along $\gamma$ defined by
\[W = U_0\begin{bmatrix}0&-1\\1&0\end{bmatrix} + U_0^\perp\begin{bmatrix}-\frac{t_1}{4}&\frac{t_1}{4}\\-\frac{t_1}{4}&-\frac{t_1}{4}\end{bmatrix}\]
which vanishes at $t_1$. This results in $\gamma$ having its first conjugate point at $t_1\approx 2.8690968494$. Since the length of $\gamma$ featuring a cut point that is no conjugate point is limited to at least $\pi$, this first conjugate point $\gamma(t_1)$ is the cut point of the geodesic $\gamma$. This example from a plane of maximal sectional curvature coincides with the observations on the injectivity radius from \Cref{sec:numexpInj}.
In summary, we have proven
\begin{theorem}
 On $St(4,2)$, for arbitrary $U_0\in St(n,p)$, the geodesic
 \[
    t\mapsto\gamma(t) = \text{Exp}_{U_0}(t\Delta), \quad \Delta = U_0\begin{bmatrix} 0&0\\0&0\end{bmatrix} + U_0^{\perp}\begin{bmatrix} \frac{1}{2}&\frac{1}{2}\\-\frac{1}{2}&\frac{1}{2} \end{bmatrix},
 \]
 with starting velocity $\Delta$ from a tangent plane section of maximal sectional curvature
 has its cut point at its first conjugate point.
 This cut point occurs at the first positive root of
 \[
    t\mapsto \left(\frac{t}{\sqrt{2}}\cos\left(\frac{t}{\sqrt{2}}\right) + \sin\left(\frac{t}{\sqrt{2}}\right)\right)
 \]
 which is at $t_1\approx 2.8690968494$.
 As this geodesic can be embedded in all Stiefel manifolds of dimensions $p\geq2, n\geq p+2$,
 the same construction gives conjugate points at the same geodesic length on all such $St(n,p)$.
\end{theorem}

The tangent vector $\Delta$ defining $\gamma$ is part of the tangent space section $\Omega = \text{span}(\Delta_1,\Delta_2)$ of maximal sectional curvature spanned by the tangent vectors
\[
    \Delta_1 = U_0 A_1 + U_0^{\perp} B_1,\,\Delta_2 = U_0 A_2 + U_0^{\perp} B_2,\;\text{ with } A_1 = A_2 = 0,\; B_1 = \frac{1}{\sqrt{2}}\begin{bmatrix}1&0\\0&1\end{bmatrix}, B_2 = \frac{1}{\sqrt{2}}\begin{bmatrix}0&1\\-1&0\end{bmatrix}.
\]
A tangent vector $\tilde\Delta\in\Omega$ of unit norm has the form 
\[\tilde\Delta(\lambda) = U_0^{\perp}\Big(\lambda B_1 + \sqrt{1-\lambda^2} B_2\Big),\quad\lambda\in\left[0,1\right].\]
By the same line of arguments as before, we obtain that all geodesics $\tilde\gamma$ defined by $U_0$ and $\tilde\Delta(\lambda)\in\Omega$ have their cut point in form of their first conjugate point at a geodesic length of $t_1\approx 2.8690968494$. 
For the five linearly independent directions $W$ that define the Jacobi fields, we choose
\begin{align*}
    W_1 &= U_0\begin{bmatrix}0&-1\\1&0\end{bmatrix},\quad W_2 = U_0^\perp\begin{bmatrix}-\lambda&\sqrt{1-\lambda^2}\\\sqrt{1-\lambda^2}&\lambda\end{bmatrix}, \quad  W_3 = U_0^\perp\begin{bmatrix}\lambda&-\sqrt{1-\lambda^2}\\\sqrt{1-\lambda^2}&\lambda\end{bmatrix},\\
    W_4 &= U_0^\perp\begin{bmatrix}\lambda&\sqrt{1-\lambda^2}\\-\sqrt{1-\lambda^2}&\lambda\end{bmatrix},\quad 
    W_5 = U_0^\perp\begin{bmatrix}\lambda&\sqrt{1-\lambda^2}\\\sqrt{1-\lambda^2}&-\lambda\end{bmatrix},
\end{align*}
for $\lambda\in(0,1)$. For $\lambda\in\{0,1\}$ the analysis can be done in the same fashion. 
We obtain linearly independent Jacobi fields. As before, the matrix function values $J_k(t)$ at any $t$ are also linearly independent as long as all $t$-dependent entries are non-zero. 
The important Jacobi fields $J_k$ for forming the Jacobi field defining the first conjugate point at $t_1$ are
\begin{align*}
    J_1(t) =(d\text{Exp}_{U_0})_{t\tilde\gamma^\prime(0)}[tW_1]=& \begin{bmatrix}
        0&-\frac{t\cos\left(\frac{t}{\sqrt{2}}\right) + \sqrt{2}\sin\left(\frac{t}{\sqrt{2}}\right)}{2}\\
        \frac{t\cos\left(\frac{t}{\sqrt{2}}\right) + \sqrt{2}\sin\left(\frac{t}{\sqrt{2}}\right)}{2} & 0\\
        \frac{\sqrt{1-\lambda^2}}{2}t\sin\left(\frac{t}{\sqrt{2}}\right)&-\frac{\lambda}{2}t\sin\left(\frac{t}{\sqrt{2}}\right)\\
        \frac{\lambda}{2}t\sin\left(\frac{t}{\sqrt{2}}\right)&\frac{\sqrt{1-\lambda^2}}{2}t\sin\left(\frac{t}{\sqrt{2}}\right)
    \end{bmatrix},\\
    J_3(t) =(d\text{Exp}_{U_0})_{t\tilde\gamma^\prime(0)}[tW_3]=& \begin{bmatrix}
        (1-2\lambda^2)t\sin\left(\frac{t}{\sqrt{2}}\right)&0\\
        0&(1-2\lambda^2)t\sin\left(\frac{t}{\sqrt{2}}\right)\\
        f(t)&g(t)\\
        -g(t)&f(t)
    \end{bmatrix},\\
    & f(t) = (2\lambda^2-1)\lambda t\cos\left(\frac{t}{\sqrt{2}}\right)+2\sqrt{2}(1-\lambda^2)\lambda\sin\left(\frac{t}{\sqrt{2}}\right),\\
    & g(t) = (2\lambda^2-1)t\sqrt{1-\lambda^2}\cos\left(\frac{t}{\sqrt{2}}\right)-2\sqrt{2}\sqrt{1-\lambda^2}\lambda^2\sin\left(\frac{t}{\sqrt{2}}\right),\\
    J_4(t) =(d\text{Exp}_{U_0})_{t\tilde\gamma^\prime(0)}[tW_4]=& \begin{bmatrix}
        -t\sin\left(\frac{t}{\sqrt{2}}\right)&0\\
        0&-t\sin\left(\frac{t}{\sqrt{2}}\right)\\
        t\lambda\cos\left(\frac{t}{\sqrt{2}}\right)&t\sqrt{1-\lambda^2} \cos\left(\frac{t}{\sqrt{2}}\right)\\
        -t\sqrt{1-\lambda^2} \cos\left(\frac{t}{\sqrt{2}}\right)&t\lambda\cos\left(\frac{t}{\sqrt{2}}\right)
    \end{bmatrix}.
\end{align*}
The Jacobi field defining the first conjugate point is given by the linear combination 
\[J(t) = J_1(t) - \frac{t_1}{4\sqrt{2}}\frac{1}{\lambda\sqrt{1-\lambda^2}}\big(J_3(t) + (1-2\lambda^2)J_4(t)\big)\]
or can equivalently be formed via the direction $W = W_1 - \frac{t_1}{4\sqrt{2}}\frac{1}{\lambda\sqrt{1-\lambda^2}}\big(W_3 + (1-2\lambda^2)W_4\big)$. Therefore, each of the geodesics have their cut point at $t_1$.

With this result, one might hope that there is a relationship of the form that geodesics defined by a tangent vector from a tangent space section with smaller sectional curvature have their cut point at a larger geodesic distance. However, this is not the case, as tangent vectors naturally do not lie in one tangent plane section only. For example, the tangent vector $\Delta$, which defines the geodesic $\gamma$ that has its cut point at $t_1$, lies not only in $\Omega$, with $\mathcal{K}(\Omega) = \frac54$, but also in the tangent plane section
\[\tilde\Omega = \Big\{a\tilde\Delta_1 + b\tilde\Delta_2\mid \tilde\Delta_1 = U_0^{\perp}\begin{bmatrix}1&0\\0&0\end{bmatrix},\,\tilde\Delta_2 = U_0^{\perp}\frac{1}{\sqrt{3}}\begin{bmatrix}0&1\\-1&1\end{bmatrix}\Big\},\]
with $\mathcal{K}(\tilde\Omega) = \frac{5}{12}$.

\begin{remark}
    There is one major difference when comparing the study of the injectivity radius on the Stiefel manifold under the canoncical metric with the Euclidean case. Under both metrics, the shortest closed geodesics have length $2\pi$. However, the sectional curvature of the Stiefel manifold under the Euclidean metric is globally bounded by $1$. Therefore, the (first) conjugate points have at least a geodesic distance of $\pi$. Hence, we avoid the analysis of conjugate points when determining the injectivity radius via Klingenberg's theorem (\Cref{thm:Klingenberg}) by stating a closed geodesic of length $2\pi$. So, the injectivity radius of the Stiefel manifold under the Euclidean metric is $\pi$. The bound on the sectional curvature was proofed in \cite[Theorem 10]{ZimmermannStoye2024}. The length of the shortest closed geodesics and thus the injectivity radius of the Stiefel manifold under the Euclidean metric was determined in \cite{zimmermannstoye2024injectivity}.
\end{remark}

\section{Summary}\label{sect:Summary}
This paper investigates the injectivity radius of the Stiefel manifold, which defines the 
size of the largest domains within which geodesics describe uniquely shortest paths everywhere on the manifold. First, we re-derive the curvature-based bound of $\sqrt{\frac{4}{5}}\pi \approx 2.81$ on the injectivity radius of the Stiefel manifold found by Rentmeesters \cite{Quentin}. The sharpness of this bound is an open question. We investigate the sharpness by investigating various random geodesics for their cut points numerically, which leads to bounding the injectivity radius by $2.87$ from above. \\
In a second part, we construct an explicit example of a cut point on the Stiefel manifold $St(4,2)$. Dimension-wise this is the first `true' Stiefel manifold in the sense that the manifolds of smaller dimensions are either isomorphic to the spheres or to the orthogonal groups 
of corresponding dimensions. Yet, it is to be expected that all extreme cases for the injectivity radius at a point, for conjugate points or for cut points already occur here, as do the extreme cases for the sectional curvature \cite{ZimmermannStoye2024}.\\
We investigate a geodesic with starting velocity from a tangent plane section of maximal sectional curvature. For this geodesic we derive a complete, in this case five-dimensional set of linearly independent Jacobi fields along the geodesic and use them to obtain the first conjugate point along the geodesic at the smallest positive root of $\left(\frac{t}{\sqrt{2}}\cos\left(\frac{t}{\sqrt{2}}\right) + \sin\left(\frac{t}{\sqrt{2}}\right)\right)$, which is at $t_1\approx 2.8690968494$. This first conjugate point is describing the cut point of the geodesic and aligns with the results from the numerical experiments. \\
Furthermore, the results are a strong support for the conjecture made by Absil and Mataigne in \cite[Conjecture 8.1]{absil2024ultimate} for the Stiefel injectivity radius in the case of the canonical metric.

\begin{gammacode}
    The Python scripts for the numerical experiments can be found in
    \begin{center}
	 	\href{https://github.com/JakobStoye/StiefelInjRadius/}{https://github.com/JakobStoye/StiefelInjRadius/}.
	\end{center}
\end{gammacode}

\begin{gammacknowledgement}
	Two days before the submitting the first preprint of this work, by coincidence, we learned about the work of \cite{absil2024ultimate} through personal communication.
    We would like to thank the authors of \cite{absil2024ultimate}, Pierre-Antoine Absil and Simon Mataigne, for a very stimulating and constructive exchange at the last minute.
\end{gammacknowledgement}




\end{document}